\documentclass[12pt]{amsart}
\usepackage{graphicx}
\usepackage{amssymb}
\usepackage{amsfonts}
\usepackage{amsmath}
\usepackage{array}
\usepackage{rotating}

\usepackage[matrix,frame]{xypic}
\newdimen\vcadre\vcadre=0.1cm 
\newdimen\hcadre\hcadre=0.1cm 
\def\GrTeXBox#1{\vbox{\vskip\vcadre\hbox{\hskip\hcadre%
      $#1$%
   \hskip\hcadre}\vskip\vcadre}}
\def\arx#1[#2]{\ifcase#1 \relax \or%
  \ar @{-}[#2]  \or%
  \ar @2{-}[#2] \or%
  \ar @{--}[#2] \or%
  \ar @2{.}[#2] \or%
  \ar @{~}[#2]  \fi}

\usepackage{tikz}
\usetikzlibrary{arrows.meta}

\def\FF{{\mathbb{F}}}
\def\Hess{{\rm Hess}}

\def\ev{{\rm ev}}

\def\<{\langle}
\def\>{\rangle}

\def\C{{\mathbb C}}

\def\NN{{\mathbb N}}

\def\KK{{\mathbb K}}

\def\SG{{\mathfrak S}}

\def\X{{\rm X}}

\def\asc{{\rm asc}}

\makeatletter
\newsavebox{\@brx}
\newcommand{\llangle}[1][]{\savebox{\@brx}{\(\m@th{#1\langle}\)}%
  \mathopen{\copy\@brx\kern-0.5\wd\@brx\usebox{\@brx}}}
\newcommand{\rrangle}[1][]{\savebox{\@brx}{\(\m@th{#1\rangle}\)}%
  \mathclose{\copy\@brx\kern-0.5\wd\@brx\usebox{\@brx}}}
\makeatother

\def\Tabvrule{\vrule width-0.4pt}       
\def\Tabhrule{\hrule \hrule height-0.4pt} 
\def\Tabstrut{\vrule height2.2ex 
                     depth0.8ex  
                     width0ex    
\relax}

\def\PasCase#1{\omit%
            $\vcenter{\hbox {\vbox to 0.4pt{}}
               \hbox{\makebox[3ex]{\Tabstrut$#1$}}}%
               \Tabvrule$}
\def\PasCasePoint{\PasCase{\cdot}}
\def\DessinCarre#1{%
    \vcenter{\hbox{}\hrule
             \hbox{\vrule\makebox[3ex]{\Tabstrut$#1$}\vrule}\Tabhrule}%
             \Tabvrule}
\def\GenRuban#1{\vcenter{\halign{&$\DessinCarre{##}$\cr#1}}\egroup}

\def\sTabvrule{\vrule width-0.4pt}
\def\sTabhrule{\hrule \hrule height-0.4pt}
\def\sTabstrut{\vrule height1.6ex depth0.6ex width0ex \relax}
\def\sDessinCarre#1{%
    \vcenter{\hbox{}\hrule
             \hbox{\vrule\makebox[2.3ex]%
                  {\sTabstrut$\scriptstyle#1$}\vrule}\sTabhrule}%
             \sTabvrule}
\def\sGenRuban#1{\vcenter{\halign{&$\sDessinCarre{##}$\cr#1}}\egroup}

\def\ruban{%
  \bgroup
  \let\ =\omit
  \let\\=\cr
  \let\x=\times
  \let\.=\PasCasePoint
  \offinterlineskip
  \GenRuban}

\def\sruban{%
  \bgroup
  \let\ =\omit
  \let\x=\times
  \let\\=\cr
  \offinterlineskip
  \sGenRuban}

\newdimen\Squaresize \Squaresize=14pt
\newdimen\Thickness \Thickness=0.5pt

\def\Square#1{\hbox{\vrule width \Thickness
   \vbox to \Squaresize{\hrule height \Thickness\vss
      \hbox to \Squaresize{\hss#1\hss}
   \vss\hrule height\Thickness}
\unskip\vrule width \Thickness}
\kern-\Thickness}

\def\Vsquare#1{\vbox{\Square{$#1$}}\kern-\Thickness}

\def\young#1{
\vbox{\smallskip\offinterlineskip
\halign{&\Vsquare{##}\cr #1}}}

\def\boxit#1#2{\setbox1=\hbox{\kern#1{#2}\kern#1}%
\dimen1=\ht1 \advance\dimen1 by #1 \dimen2=\dp1 \advance\dimen2 by #1
\setbox1=\hbox{\vrule height\dimen1 depth\dimen2\box1\vrule}%
\setbox1=\vbox{\hrule\box1\hrule}%
\advance\dimen1 by .4pt \ht1=\dimen1
\advance\dimen2 by .4pt \dp1=\dimen2 \box1\relax}

\def\PC{{\rm PC\,}}

\def\taille{.5}

\title[]{Acyclic orientations and Hessenberg varieties}

\author[ J.-C.~Novelli and J.-Y.~Thibon]%
{Jean-Christophe Novelli and Jean-Yves Thibon}

\address[] {LIGM, Universit\'e
Gustave-Eiffel, CNRS, ENPC, ESIEE-Paris \\
5 Boulevard Descartes \\Champs-sur-Marne \\77454 Marne-la-Vall\'ee cedex 2 \\
FRANCE}
\email[Jean-Christophe Novelli]{jean-christophe.novelli@univ-eiffel.fr}
\email[Jean-Yves Thibon]{jean-yves.thibon@univ-eiffel.fr} 


\date{}

\begin{document}

\begin{abstract}
	We exhibit a bijection between acyclic orientations of a Dyck graph and
	Tymoczko cells of a regular nilpotent Hessenberg variety. This implies
	the Shareshian-Wachs formula for the sum of the coefficients of the chromatic
	quasi-symmetric function of a Dyck graph in the elementary basis.
\end{abstract}

\maketitle

\section{Introduction}

In 1995, Stanley introduced a symmetric analogue $X_G$ of the chromatic polynomial
of a graph $G$, and proved among other things that the sum of its coefficients
in the elementary basis is equal to the number of acyclic orientations of the graph \cite{St}.

Some twenty years later, Shareshian and Wachs introduced a $t$-analogue $X_G(t)$ of Stanley's
chromatic symmetric function, which is in general only a quasi-symmetric function \cite{SW}.
However, for a very important family of graphs, the incomparability graphs of  natural unit
interval orders, called here for short Dyck graphs, this quasi-symmetric function turns
out to be symmetric, and Shareshian-Wachs proved that the sum of its coefficients on
the elementary basis counts acyclic orientations of $G$ by the number of ascending edges.

Actually, Stanley conjectured that for Dyck graphs, $X_G$ is $e$-positive,
and Shareshian-Wachs conjectured the stronger statement that $X_G(t)$ is also $e$-positive.

Another conjecture of Shareshian-Wachs, now a theorem \cite{BrCh}, relates $X_G(t)$ to the 
geometry of Hessenberg varieties: $\omega X_G(t^2)$ is the equivariant Poincaré polynomial
of a regular semisimple Hessenberg variety associated with the graph $G$. The $e$-positivity
conjecture means then  that the representation of the symmetric group on the equivariant cohomology
is actually a permutation representation.

But $\omega X_G(t)$ also encodes the Poincaré polynomials of nilpotent Hessenberg varieties,
which are obtained by taking its scalar product with  certain Hall-Littlewood functions.
In the so-called regular case, these Hall-Littlewood functions reduce to the  $h_n$, so that the Poincaré
polynomial is then precisely the sum of the coefficients of $X_G(t)$ on the elementary basis.

Tymoczko has given cell decompositions (affine pavings) of various types of Hessenberg
varieties \cite{Tym} allowing a purely combinatorial calculation of their Poincaré polynomials.

The aim of this note is to give a geometric proof of the Shareshian-Wachs formula by exhibiting
a bijection between Tymoczko cells of a regular nipotent Hessenberg variety and acyclic orientations
of the associated Dyck graph.

{\bf Acknowlegements. } This research has been partially supported by the project CARPLO
of the Agence Nationale de la recherche (ANR-20-CE40-0007).


\section{Background}
\subsection{Chromatic quasi-symmetric functions}
Let $G$ be a simple undirected graph on the set of vertices $V(G)=\{1,\ldots,n\}$, 
and let $E(G)$ denote  the set of edges of $G$.
A coloring of $G$ is a map $c:\ V(G)\rightarrow \NN^*$, which can be
identified with a word $c_1c_2\cdots c_n$.
A coloring is proper if $c_i\not=c_j$ whenever $\{i,j\}\in E(G)$. We denote by
$C(G)$ the set of proper colorings of $G$.

The chromatic quasi-symmetric function of $G$
is defined by \cite{SW}
\begin{equation}
\X_G(t,X) = \sum_{c\in C(G)}t^{\asc_G(c)}x_{c_1}x_{c_2}\cdots x_{c_n}
        =  \sum_{c\in \PC(G)}t^{\asc_G(c)}M_{\ev(c)}(X),
\end{equation}
where $\PC(G)$ denotes the set of proper {\it packed colorings}
({\it i.e.,} whose colors form an interval $[1,r]$), $\asc_G(c)$ is the
number of edges $\{i<j\}$ such that $c_i<c_j$, and $\ev(c)$ is the evaluation (or content)
of $c$, that is, the integer composition recording the number of occurences of each value of $c$.

\subsection{Dyck graphs}
Dyck graphs $G$ are simple undirected graphs with vertices labelled $1,\ldots,n$,
characterized by the property that if there is an edge $\{i,j\}$ with $i<j$,
then all the $\{i',j'\}$ with $i\leq i'<j'\leq j$ are also edges of $G$.
The number of such graphs is the Catalan number $c_n$. These are the
incomparability graphs of certain posets $P$, known as unit interval orders \cite{SW}.

Dyck graphs on $n$ vertices are in bijection with Dyck paths of length $2n$, whence their name.

\subsection{Hessenberg varieties}
The Hessenberg variety $\Hess(X,h)$ associated with a function
$h:\ [n]\rightarrow [n]$ and a matrix $X\in M_n(\KK)$ is the set of complete flags
$V_1\subset V_2\subset\cdots\subset V_n$ of $\KK^n$ such that $XV_i \subseteq V_{h(i)}$.

The most studied cases are:

a) The Hessenberg function is increasing and satisfies $h(i)\ge i$, which encompasses

- the semisimple case: $X=S$, a diagonal matrix with distinct eigenvalues,

- the regular case: the Jordan blocks of $X$ have distinct eigenvalues,

- the  nilpotent case: $X=N_\mu$, a nilpotent matrix whose Jordans blocks are specified by a partition $\mu$.
When $\mu=(n)$, one speaks of a regular nilpotent variety.

b) The  Hessenberg function is increasing and satisfies $h(i)< i$.

The only studied case is for $X=N_\mu$, the so-called ad-nilpotent case..

Hessenberg functions of type a) are in bijection with Dyck graphs: the edges are $\{i,j\}$
for $i<j\le h(i)$.

$h(i)-i$ is the number of edges $(i<j)$.

Hessenberg functions $h$ of type a)  also correspond to partitions $\mu$  contained in the staircase $(n-1,\ldots,2,1)$ and $312$-avoiding permutations $w$.
The partition $\mu$ is the complement of $h$ in an $n\times n$ square diagram (i.e., $\mu_i = n-h(i)$).
The edges of the corresponding graph $G$ are the coordinates of the empty cells above the diagonal, and the code of $w$ is $c_i=h(i)-i$.

\medskip
{\footnotesize
For example, $h=23555$ corresponds to $\mu=221$, the edges of $G$ are
$(1,2),(2,3),(3,4),(3,5),(4,5)$, and  the code of $w$ is $c=11210$ so that $w=23541$:
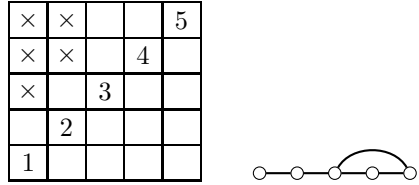
\begin{equation*}
\young{\times &\times & & & 5\cr
       \times & \times & &4& \cr
       \times & & 3 &&\cr
        & 2 &&&\cr
        1&&&&\cr}
\quad\quad
\begin{tikzpicture}
\begin{scope}[every node/.style={circle,scale=.5,fill=white,draw}]
    \node (A) at (0,0) {};
    \node (B) at (1*\taille,0) {};
    \node (C) at (2*\taille,0) {};
    \node (D) at (3*\taille,0) {};
    \node (E) at (4*\taille,0) {};
\end{scope}

\begin{scope}[>={Stealth[black]},
              every edge/.style={draw=black,thick}]
    \path [-] (A) edge (B);
    \path [-] (B) edge (C);
    \path [-] (C) edge (D);
    \path [-] (D) edge (E);
    \path [-] (C) edge[bend left=60] (E);
\end{scope}
\end{tikzpicture}
\end{equation*}
}

\subsection{Affine pavings and Poincar\'e polynomials}

The varieties of the cases a) admit affine pavings, i.e., decompositions as disjoint unions of affine spaces (cells),
parametrized by Tymoczko's configurations \cite{Tym}. 
These are standard fillings of Ferrers diagrams\footnote{We take the French convention for Ferrers diagrams.}, on which one can read the dimensions of the cells 
as a number of special inversions, depending on the function $h$. 

These configurations give therefore the number of $\FF_q$-rational points of the variety as well as its Poincar\'e polynomial for $\KK=\C$.

The Tymoczko configurations associated with a Hessenberg function $h$ and a nilpotent matrix of type $\mu$
are the standard fillings of the Ferrers diagram of $\mu$ in which an adjacency $ij$ in a row is allowed only if 
$i\le h(j)$. An $h$-inversion is a pair $(a,b)$ with  $a<b$ and either  $b$ is above $a$ in its column, either $b$ is strictly to the left of $a$. 
One demands moreover that if $a$ is immediately followed by $c$ in its row, $b\le h(c)$.

The dimension of a cell is the number of $h$-inversions of its configuration.

\medskip
{\footnotesize
For example, with $h=233$ we have three nilpotent varieties, whose affine cells with their dimensions are as follows.

Nilpotent $(3)$:
$$
{\young{1&2&3\cr}\atop 0}\quad
{\young{1&3&2\cr}\atop 1}\quad
{\young{2&1&3\cr}\atop 1}\quad
{\young{3&2&1\cr}\atop 2}
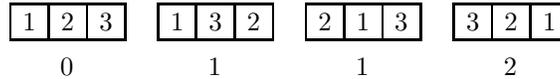
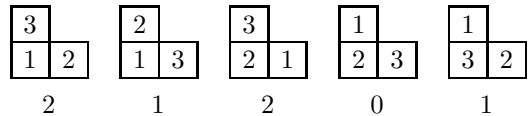
$$

Nilpotent $(2,1)$:
$$
{\young{3\cr 1&2\cr} \atop 2}\quad
{\young{2\cr 1&3\cr} \atop 1}\quad
{\young{3\cr 2&1\cr} \atop 2}\quad
{\young{1\cr 2&3\cr} \atop 0}\quad
{\young{1\cr 3&2\cr} \atop 1}
$$

Nilpotent $(1,1,1)$:
$$
{\young{3\cr 2\cr 1\cr}\atop 3 }\quad
{\young{2\cr 3\cr 1\cr}\atop 2 }\quad
{\young{3\cr 1\cr 2\cr}\atop 2 }\quad
{\young{1\cr 3\cr 2\cr}\atop 1 }\quad
{\young{2\cr 1\cr 3\cr}\atop 1 }\quad
{\young{1\cr 2\cr 3\cr}\atop 0 }
$$

In this very simple case, one can describe the variety explicitly, and compute its number of $\FF_q$-points.
We want to calculate the number $d_\mu(q)$ of complete flags of
 $\mathbb{F}_q^3$ such that  $NV_i\subseteq V_{h(i)}$ for a nilpotent matrix  $N$ whose Jordan form is specified by a partition  $\mu$ of 3.
The only condition is thus $NV_1\subseteq V_2$.

If $N=N_{111}=0$, the condition is void, and thus satisfied by all flags.
Therefore, $d_{111}(q)=[3]_q!$.

If $\mu=(21)$, set $v=x_1e_1+x_2e_2+x_3e_3$ and $V_1=\langle v\rangle$.
Then, $Nv=x_2e_1\in V_2$ and either  $x_2=0$ so that $V_1$ is an arbitrary line
of the plane $\langle e_1,e_3\rangle$, which gives $(q+1)$ possibilities, and $V_2$ is an arbitrary plane containing $V_1$, hence again $(q+1)$ choices, 
which yields a total of
$(q+1)^2$ for the case $x_2=0$.

If $x_2\not=0$, then $V_2=\langle e_1,e_2+\alpha e_3\rangle$, with $q$ choices for $\alpha$,
and $v=xe_1+y(e_2+\alpha e_3)$ with $y\not=0$ hence $q(q-1)$ choices for $(x,y)$, to be divided by
$q-1$, hence $q$ choices of $V_1$ for each choice of $V_2$, in total $q^2$ choices for the case
$x_2\not =0$.
The Poincar\'e polynomial is then $d_{21}(q)=(1+q)^2+q^2=1+2q+2q^2$ for this case.

If $\mu=(3)$, either $Nv=0$ in which case $V_1=\ker N=\langle e_1\rangle$ and $V_2$ is an arbitrary plane containing $V_1$, which gives $q+1$ choices, either 
$Nv=x_2e_1+x_3e_2\not=0$, in which case one has to choose
$(x_1,x_2,x_3)\not=(x_1,0,0)$, whence $q$ choices of $x_1$
and $q^2-1$ choices of  $(x_2,x_3)$, to be divided by  $q-1$, which yields $q^2+q$.
The Poincar\'e polynomial is therefore $d_{3}(q)=q+1+q^2+q=1+2q+q^2$.

We have recovered the values predicted by Tymoczko's configurations.

By an unpublished result of McPherson and Tymoczko \cite{OWR-Tym}, later extended to other types by
Ji and Precup \cite{JP},
they are also
given by the scalar products with the modified Hall-Littlewood functions
\begin{equation}
d_\mu(q)=\langle\omega X_G,\tilde Q'_\mu\rangle\quad\text{where}\ \tilde Q'_\mu=\sum_\lambda \tilde K_{\lambda\mu}(q)s_\lambda
\end{equation}
Here,
\begin{equation}
\omega X_G = qh_{21}+(1+q+q^2)h_3, 
\end{equation}
$$\tilde Q'_3 = m_{111}+m_{21}+m_3,\quad \tilde Q'_{21}=(2q+1)m_{111}+(q+1)m_{21}+m_3,\quad \tilde Q'_{111}=[3]_q m_{111}+[3]_q m_{21}+m_3$$
and we can see that the scalar products give back the expected results.

}

\section{Acyclic orientations and the regular nilpotent case}

For  $\mu=(n)$, according to the above mentioned result of McPherson-Tymoczko, the Poincar\'e
polynomial is
\begin{equation}
P_{h,(n)}(t)=\langle\omega X_G(t),\tilde Q'_n\rangle
 =\langle\omega X_G(t),h_n\rangle
=\sum_{\lambda\vdash n}\langle\omega X_G(t),m_\lambda\rangle
\end{equation}
since $\tilde Q'_n=h_n$. 

It is therefore also the sum of the coefficients of
$X_G(t)$ in the $e$-basis, which, according to Shareshian-Wachs, is equal to the generating
polynomial of acyclic orientations of $G$ counted by the number of increasing arcs.

We shall now reprove this result, as well as the special case $\mu=(n)$ of the McPherson-Tymoczko formula,
by establishing a statistic-preserving  bijection between acyclic orientations and Tymoczko cells.

Each acyclic orientation $o$ determines a partial order $P(o)$ on $[n]$. We shall show that each order $P(o)$
has a unique minimal linear extension $\sigma_o$, which can be read as a Tymoczko cell of shape $(n)$, and whose
dimension is equal to the number of increasing arcs of $o$.

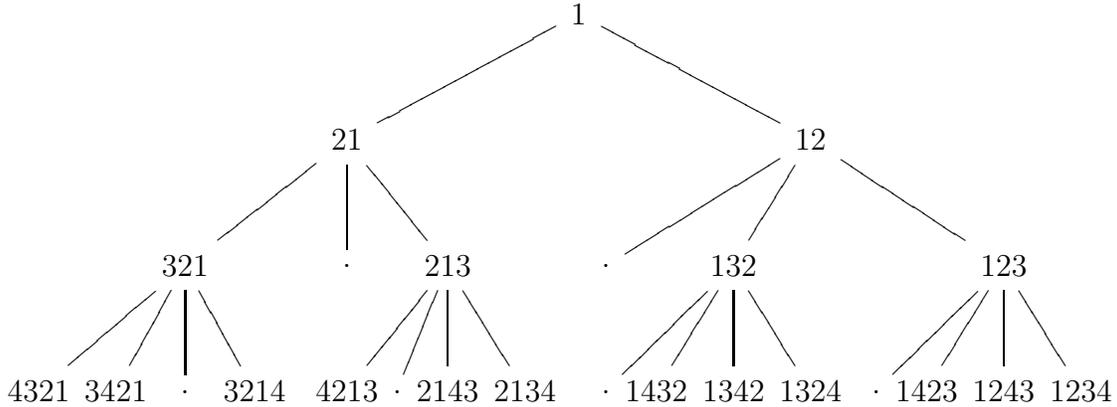
\begin{figure}[!t]
	\centering
\newdimen\vcadre\vcadre=0.2cm 
\newdimen\hcadre\hcadre=0.2cm 
$\xymatrix@R=1.0cm@C=-2mm{
&&&&&&&&& *{\GrTeXBox{1}}\arx1[lllld]\arx1[rrrrd]& \\
&&&&& *{\GrTeXBox{21}} \arx1[llld]\arx1[d]\arx1[rrd]
&&&&&&&& *{\GrTeXBox{12}} \arx1[llld]\arx1[ld]\arx1[rrrrd]&
\\
&&  *{\GrTeXBox{321}} \arx1[lld]\arx1[ld]\arx1[d]\arx1[rd]
&&& *{\GrTeXBox{.}}
&& *{\GrTeXBox{213}} \arx1[lld]\arx1[ld]\arx1[d]\arx1[rd]
&&& *{\GrTeXBox{.}}
&& *{\GrTeXBox{132}} \arx1[lld]\arx1[ld]\arx1[d]\arx1[rd]
&&&&& *{\GrTeXBox{123}} \arx1[lld]\arx1[ld]\arx1[d]\arx1[rd]
\\
*{\GrTeXBox{4321}} &
*{\GrTeXBox{3421}} &
*{\GrTeXBox{.}} &
*{\GrTeXBox{3214}} &
*{\GrTeXBox{}} &
*{\GrTeXBox{4213}} &
*{\GrTeXBox{.}} &
*{\GrTeXBox{2143}} &
*{\GrTeXBox{2134}} &
*{\GrTeXBox{}} &
*{\GrTeXBox{.}} &
*{\GrTeXBox{1432}} &
*{\GrTeXBox{1342}} &
*{\GrTeXBox{1324}} &
*{\GrTeXBox{}} &
*{\GrTeXBox{.}} &
*{\GrTeXBox{1423}} &
*{\GrTeXBox{1243}} &
*{\GrTeXBox{1234}}
\\
}$

	\caption{Generation tree for $h=2444$.}\label{fig1}
\end{figure}

This property is proved by examining the generation tree of Tymoczko cells of shape $(n)$ (see Figure \ref{fig1}).
One puts 1 at the root, and the descendants of each node correspond to the different possibilities
of inserting $k$ in a permutation of $[k-1]$.

On can insert $k>i$ to the left of $i$ only if $k\le h(i)$, which amounts to the condition
 $\{i,k\}\not\in E(G)$. In this case, the insertions
$k'i$ will also be forbidden for all $k'>k$, by the definition of Dyck graphs.

Thus, each permutation $\sigma$  of a same level  $k-1$ 
has the same number of descendants, since one can only insert $k$ at the end or before an $i$ to which is is connected by an edge.
These insertions do not alter the $h$-inversions which were already present in $\sigma$, and add successively, starting from the right,
$0,1,2,\ldots b_k$ inversions between
 $k$ and the predecessors of those to which it is connected by an edge, whose number is
$b_k=h(k)-k$.

The generating polynomial of $h$-inversions is therefore
\begin{equation}
P_{h,(n)}(t)=\prod_{i=1}^n [h(i)-i+1]_t.
\end{equation}

On another hand, the acyclic orientations of $G$ determine a partition of
$\mathfrak{S}_n$ whose blocks consist of the linear extensions of the $P(o)$ (Figure \ref{fig2}). 

Two permutations $\sigma=uijv$ and $\tau=ujiv$ are in the same class iff
$\{i,j\}\not\in E(G)$. 
Indeed, is there were a path from
$i$ to $j$ in the orientation $o$, this path would have at least length 2, and $i,j$ could not be adjacent in a linear extension. 

The local minima of the classes coincide with the Tymoczko configurations of shape $(n)$. Two such configurations
cannot be linear extensions of the same $P(o)$: if 
$\alpha\in\mathfrak{S}_{k-1}$ is their first common ancestor in the generation tree,
the descendants of $\alpha$ are distinguished by the relative positions of
$k$ and the  $i<k$ such that $\{i,k\}\in E(G)$, 
which amounts to changing the orientation of one of these edges.

There is therefore a unique minimum in each class, which determines a bijection between acyclic orientations $o$
and configurations $\tau$, sending the number of increasing arcs of $o$
to the number of $h$-inversions of $\tau$.

Thus, if one starts with the Tymoczko decomposition, one obtains that the Poincar\'e polynomial is given by
$P_{h,(n)}=\langle\omega X_G(t),\tilde Q'_n\rangle$, 
that is, it is given by the product of $t$-integers above, and that it counts acyclic orientations by increasing arcs.

{\footnotesize
\begin{figure}[!t]
	\centering
\begin{equation*}
\begin{split}
\begin{tikzpicture}
\begin{scope}[every node/.style={circle,scale=.5,fill=white,draw}]
   \node (A) at (0,0) {}; \node (B) at (1*\taille,0) {};
    \node (C) at (2*\taille,0) {}; \node (D) at (3*\taille,0) {};
\end{scope}
\begin{scope}[>={Stealth[black]}, every edge/.style={draw=black,thick}]
    \path [->] (A) edge (B); \path [->] (B) edge (C);
    \path [->] (C) edge (D);
    \path [->] (B) edge[bend left=60] (D);
\end{scope}
\end{tikzpicture}
	& ~\hskip.1cm: {\bf 1234} \\
\begin{tikzpicture}
\begin{scope}[every node/.style={circle,scale=.5,fill=white,draw}]
   \node (A) at (0,0) {}; \node (B) at (1*\taille,0) {};
    \node (C) at (2*\taille,0) {}; \node (D) at (3*\taille,0) {};
\end{scope}
\begin{scope}[>={Stealth[black]}, every edge/.style={draw=black,thick}]
    \path [->] (A) edge (B); \path [->] (B) edge (C); \path [<-] (C) edge (D);
    \path [->] (B) edge[bend left=60] (D);
\end{scope}
\end{tikzpicture}
	& ~\hskip.1cm: {\bf 1243} \\
\begin{tikzpicture}
\begin{scope}[every node/.style={circle,scale=.5,fill=white,draw}]
   \node (A) at (0,0) {}; \node (B) at (1*\taille,0) {};
    \node (C) at (2*\taille,0) {}; \node (D) at (3*\taille,0) {};
\end{scope}
\begin{scope}[>={Stealth[black]}, every edge/.style={draw=black,thick}]
    \path [->] (A) edge (B); \path [->] (B) edge (C); \path [<-] (C) edge (D);
    \path [<-] (B) edge[bend left=60] (D);
\end{scope}
\end{tikzpicture}
	& ~\hskip.1cm: {\bf 1423}, 4123 \\
\begin{tikzpicture}
\begin{scope}[every node/.style={circle,scale=.5,fill=white,draw}]
   \node (A) at (0,0) {}; \node (B) at (1*\taille,0) {};
    \node (C) at (2*\taille,0) {}; \node (D) at (3*\taille,0) {};
\end{scope}
\begin{scope}[>={Stealth[black]}, every edge/.style={draw=black,thick}]
    \path [->] (A) edge (B); \path [<-] (B) edge (C); \path [->] (C) edge (D);
    \path [->] (B) edge[bend left=60] (D);
\end{scope}
\end{tikzpicture}
	& ~\hskip.1cm: {\bf 1324}, 3124 \\
\begin{tikzpicture}
\begin{scope}[every node/.style={circle,scale=.5,fill=white,draw}]
   \node (A) at (0,0) {}; \node (B) at (1*\taille,0) {};
    \node (C) at (2*\taille,0) {}; \node (D) at (3*\taille,0) {};
\end{scope}
\begin{scope}[>={Stealth[black]}, every edge/.style={draw=black,thick}]
    \path [->] (A) edge (B); \path [<-] (B) edge (C); \path [->] (C) edge (D);
    \path [<-] (B) edge[bend left=60] (D);
\end{scope}
\end{tikzpicture}
	& ~\hskip.1cm: {\bf 1342}, 3142, 3412 \\
\begin{tikzpicture}
\begin{scope}[every node/.style={circle,scale=.5,fill=white,draw}]
   \node (A) at (0,0) {}; \node (B) at (1*\taille,0) {};
    \node (C) at (2*\taille,0) {}; \node (D) at (3*\taille,0) {};
\end{scope}
\begin{scope}[>={Stealth[black]}, every edge/.style={draw=black,thick}]
    \path [->] (A) edge (B); \path [<-] (B) edge (C); \path [<-] (C) edge (D);
    \path [<-] (B) edge[bend left=60] (D);
\end{scope}
\end{tikzpicture}
& ~\hskip.1cm: {\bf 1432}, 4132, 4312 \\
\begin{tikzpicture}
\begin{scope}[every node/.style={circle,scale=.5,fill=white,draw}]
   \node (A) at (0,0) {}; \node (B) at (1*\taille,0) {};
    \node (C) at (2*\taille,0) {}; \node (D) at (3*\taille,0) {};
\end{scope}
\begin{scope}[>={Stealth[black]}, every edge/.style={draw=black,thick}]
    \path [<-] (A) edge (B); \path [->] (B) edge (C);
    \path [->] (C) edge (D);
    \path [->] (B) edge[bend left=60] (D);
\end{scope}
\end{tikzpicture}
& ~\hskip.1cm: {\bf 2134}, 2314, 2341 \\
\begin{tikzpicture}
\begin{scope}[every node/.style={circle,scale=.5,fill=white,draw}]
   \node (A) at (0,0) {}; \node (B) at (1*\taille,0) {};
    \node (C) at (2*\taille,0) {}; \node (D) at (3*\taille,0) {};
\end{scope}
\begin{scope}[>={Stealth[black]}, every edge/.style={draw=black,thick}]
    \path [<-] (A) edge (B); \path [->] (B) edge (C); \path [<-] (C) edge (D);
    \path [->] (B) edge[bend left=60] (D);
\end{scope}
\end{tikzpicture}
& ~\hskip.1cm: {\bf 2143}, 2413, 2431 \\
\begin{tikzpicture}
\begin{scope}[every node/.style={circle,scale=.5,fill=white,draw}]
   \node (A) at (0,0) {}; \node (B) at (1*\taille,0) {};
    \node (C) at (2*\taille,0) {}; \node (D) at (3*\taille,0) {};
\end{scope}
\begin{scope}[>={Stealth[black]}, every edge/.style={draw=black,thick}]
    \path [<-] (A) edge (B); \path [->] (B) edge (C); \path [<-] (C) edge (D);
    \path [<-] (B) edge[bend left=60] (D);
\end{scope}
\end{tikzpicture}
& ~\hskip.1cm: {\bf 4213}, 4231 \\
\begin{tikzpicture}
\begin{scope}[every node/.style={circle,scale=.5,fill=white,draw}]
   \node (A) at (0,0) {}; \node (B) at (1*\taille,0) {};
    \node (C) at (2*\taille,0) {}; \node (D) at (3*\taille,0) {};
\end{scope}
\begin{scope}[>={Stealth[black]}, every edge/.style={draw=black,thick}]
    \path [<-] (A) edge (B); \path [<-] (B) edge (C); \path [->] (C) edge (D);
    \path [->] (B) edge[bend left=60] (D);
\end{scope}
\end{tikzpicture}
& ~\hskip.1cm: {\bf 3214}, 3241 \\
\begin{tikzpicture}
\begin{scope}[every node/.style={circle,scale=.5,fill=white,draw}]
   \node (A) at (0,0) {}; \node (B) at (1*\taille,0) {};
    \node (C) at (2*\taille,0) {}; \node (D) at (3*\taille,0) {};
\end{scope}
\begin{scope}[>={Stealth[black]}, every edge/.style={draw=black,thick}]
    \path [<-] (A) edge (B); \path [<-] (B) edge (C); \path [->] (C) edge (D);
    \path [<-] (B) edge[bend left=60] (D);
\end{scope}
\end{tikzpicture}
& ~\hskip.1cm: {\bf 3421} \\
\begin{tikzpicture}
\begin{scope}[every node/.style={circle,scale=.5,fill=white,draw}]
   \node (A) at (0,0) {}; \node (B) at (1*\taille,0) {};
    \node (C) at (2*\taille,0) {}; \node (D) at (3*\taille,0) {};
\end{scope}
\begin{scope}[>={Stealth[black]}, every edge/.style={draw=black,thick}]
    \path [<-] (A) edge (B); \path [<-] (B) edge (C); \path [<-] (C) edge (D);
    \path [<-] (B) edge[bend left=60] (D);
\end{scope}
\end{tikzpicture}
& ~\hskip.1cm: {\bf 4321} \\
\end{split}
\end{equation*}

	\caption{Partition of $\SG_4$ for $h=2444$}\label{fig2}
\end{figure}
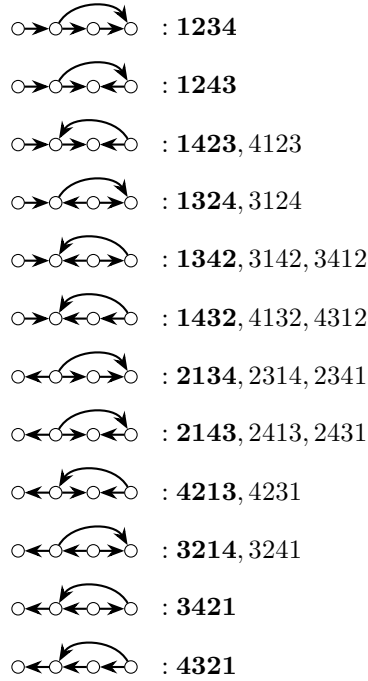
}

This yields 1) a proof of the McPherson-Tymoczko formula for the regular nilpotent case, 2) an elementary proof of the formula for the Poincar\'e polynomial
in this case, 3) a geometric or combinatorial proof of the Shareshian-Wachs formula for the sum of the coefficients on the $e$-basis



\footnotesize
\begin{thebibliography}{aa}



%
\bibitem{BrCh} P. Brosnan and T. Y. Chow, Unit interval orders and the dot action on the cohomology of regular semisimple Hessenberg varieties,
Advances in Math.  329 (2018),  955--1001
%
%
		
%
\bibitem{JP}C. Ji, M. Precup, 
Hessenberg varieties associated to ad-nilpotent ideals, Comm. Algebra 50 (2022), no. 4, 1728--1749
%



%
\bibitem{Mcd} I.G. Macdonald,  Symmetric functions and Hall
polynomials, 2nd ed., Oxford University Press, Oxford, 1995.
%
%

\bibitem{SW} J. Shareshian, M. Wachs,
Chromatic quasisymmetric functions, Advances in Math. 295 (2016), 497--551.
%
\bibitem{St} R.P. Stanley, A symmetric function generalization of the chromatic
polynomial of a graph, Adv. Math. 111 (1995), no. 1, 166--194.

\bibitem{Tym} J. Tymoczko, Paving Hessenberg varieties by affines, Selecta Math. (N.S.) 13 (2007), 353--367.

\bibitem{OWR-Tym} J. Tymoczko, An introduction to Hessenberg varieties,
	Oberwolfach Report Report No. 46/2019 (DOI: 10.4171/OWR/2019/46), 2880--2883.
%
\end{thebibliography}
\end{document}